\documentclass{amsart}

\newcommand{\ignore}[1]{}

\usepackage{graphicx}

\newcommand{\dfn}[1]{\textbf{#1}}

\newtheorem{theorem}{Theorem}

\newtheorem{corollary}[theorem]{Corollary}

\newtheorem{lemma}[theorem]{Lemma}
\newcommand{\lem}[1]{ \begin{lemma} {#1} \end{lemma} }
\newtheorem{proposition}[theorem]{Proposition}

\newtheorem{definition}{Definition}

\newtheorem{conjecture}{Conjecture}

\newtheorem{question}{Question}

\newtheorem*{NOTE}{Note}

\newtheorem*{REMARK}{Remark}

\newtheorem*{LEMMA}{Lemma}
\newcommand{\lemx}[1]{ \begin{LEMMA} {#1} \end{LEMMA} }
\newtheorem*{COROLLARY}{Corollary}
\newcommand{\corx}[1]{ \begin{COROLLARY} {#1} \end{COROLLARY} }
\newtheorem*{PROPOSITION}{Proposition}
\newcommand{\propx}[1]{ \begin{PROPOSITION} {#1} \end{PROPOSITION} }

\newcommand{\pf}[1]{\begin{proof}{#1}\end{proof}}

\newcommand{\figref}[1]{Figure~\ref{#1}}
\newcommand{\bibtitle}[1]{\emph{#1}}

\newcommand{\rr}{{\mathbb{R}}}
\newcommand{\zz}{{\mathbb{Z}}}

\newcommand{\lk}{\mathrm{lk}}

\newcommand{\qd}{(C_1, C_2, C_3, C_4)}
\newcommand{\zzt}{in $\zz$ or $\zz_2$}
\newcommand{\dl}{double-linked}
\newcommand{\dless}{$D_4$-less embedding}
\newcommand{\ty}{\nabla\mathrm{Y}}

\begin{document}{

\title{An Algorithm for Detecting Intrinsically Knotted Graphs}
\author {Jonathan Miller}
\address{
Department of Mathematics, 
Occidental College,
Los Angeles, CA 90041, USA.
}
\author{Ramin Naimi}
\address{
Department of Mathematics, 
Occidental College,
Los Angeles, CA 90041, USA.
}
\date{3 July 2013}
\thanks{The second author was supported in part by NSF grant DMS-0905300.}

\begin{abstract}
We describe an algorithm that recognizes
some (possibly all) intrinsically knotted (IK) graphs,
and can help find knotless embeddings 
for graphs that are not IK.
The algorithm, implemented as a Mathematica program,
has already been used in \cite{gmn} 
to greatly expand the list of known  minor minimal IK graphs,
and to find knotless embeddings for some graphs 
that had previously resisted attempts to classify them as IK or non-IK.
\end{abstract}

\maketitle

\section{Introduction}

All graphs in the paper are finite.
A graph is \dfn{intrinsically knotted} (IK)
if every embedding of it in $S^3$ contains a nontrivial knot.
In the early 1980s Conway and Gordon~\cite{cg}
showed that $K_7$, the complete graph on seven vertices, is IK.
A graph $H$ is a \dfn{minor} of another graph $G$
if $H$ can be obtained from a subgraph of $G$
by contracting some edges.
A graph $G$ with a given property
is said to be \dfn{minor minimal} with respect to that property
if no proper minor of $G$ has that property.
Robertson and Seymour's Graph Minor Theorem~\cite{rs}
says that in any infinite set of graphs,
at least one is a minor of another.
It follows that for any property whatsoever,
there are only finitely many graphs that are minor minimal with respect to that property.
In particular, there are only finitely many minor minimal IK graphs.
The set of all minor minimal IK graphs is also known as 
the \dfn{obstruction} set of graphs for the knotless embedding property
(a graph embedding in $S^3$ is knotless if it has no nontrivial knots).
It is well known and easy to show that $G$ is IK iff $G$ contains a minor minimal IK graph as a minor.
Hence, if we had a list of all minor minimal IK graphs,
then we could determine whether or not any given graph $G$ is IK
by testing whether $G$ contains one of the graphs in this list as a minor.
But, for about the past three decades, finding such a list has remained an open problem.
Our algorithm, which does not rely on having such a list available, 
can recognize some (possibly all) IK graphs.

Prior to this algorithm, finding new minor minimal IK graphs was quite difficult.
Initially, the only known minor minimal IK graphs were
$K_7$ and the thirteen graphs obtained from it by $\ty$ moves (\cite{cg}, \cite{ks}).
About two decades later the list expanded to include 
$K_{3,3,1,1}$ and the 25 graphs obtained from it by $\ty$ moves (\cite{foisy-D4}, \cite{ks}).
And a few years later Foisy \cite{foisy-oneMoreGraph} discovered one more minor minimal IK graph.
These 41 graphs were the only minor minimal IK graphs known before this algorithm.
With the aid of our algorithm, Goldberg, Mattman, and Naimi \cite{gmn}
have found over 1600 new IK graphs, 
most of which seem likely to be minor minimal,
and over 200 of which have been proven to be so.
And, judging from the ease with which these new graphs were found using the computer program,
 it seems there may be lots more minor minimal IK graphs.
Thus, the obstruction set for the knotlless embedding property
is much larger than it seemed (or was hoped) it might be a decade ago.
This will hopefully give researchers some direction in
what approaches to try (or, at least, not to try)
in trying to classify IK graphs.
Our algorithm was also used in \cite{gmn} in finding knotless embeddings
for some graphs that had previously resisted attempts to determine whether or not they are IK.

Our algorithm is not polynomial time
since it first needs to find all the cycles in the given graph,
and the number of cycles can grow exponentially with the number of edges.
In fact, a good predictor of computation time for our computer program 
seems to be the number of edges in the graph.
For graphs with 28 or fewer edges, 
the program usually took from about a few seconds to a few hours to run 
on a typical personal computer in 2009
(and a large portion of the time was usually spent on just finding the cycles).
More specifically, the computer program was run 
on at least the following graphs listed in \cite{gmn} and  \cite{gmn-appendix}:
all graphs in the $K_7$, $K_{3,3,1,1}$, and $E_9 + e$ families  (188, all together);
the 25 parentless (Y-less) graphs and a few dozen other graphs in the $G_{9,28}$ family;
and the graph $G_{14,25}$ plus a few more in its family.
We found that for some  graphs with about 30 or more edges, 
the program could take days (or much longer) to run.

Before describing our algorithm,
we need to introduce some more notation and terminology.
The linking number of two disjoint cycles $C_1$ and $C_2$
in a spatial graph (i.e., a graph embedded in $S^3$)
is denoted by $\lk(C_1, C_2)$. 
If $K$ is a knot in $S^3$, 
$a_2(K)$ denotes the second coefficient of the Conway polynomial of $K$.
The graph depicted in \figref{fig-D4} is denoted by $D_4$.
It has four vertices and eight edges.
The cycles $C_1, \cdots, C_4$ 
are always assumed to be labeled in the order shown, 
so that
$C_1 \cap C_3 = \emptyset$ and 
$C_2 \cap C_4 = \emptyset$.

\begin{figure}[ht]
 \centering
 \includegraphics[width=40mm]{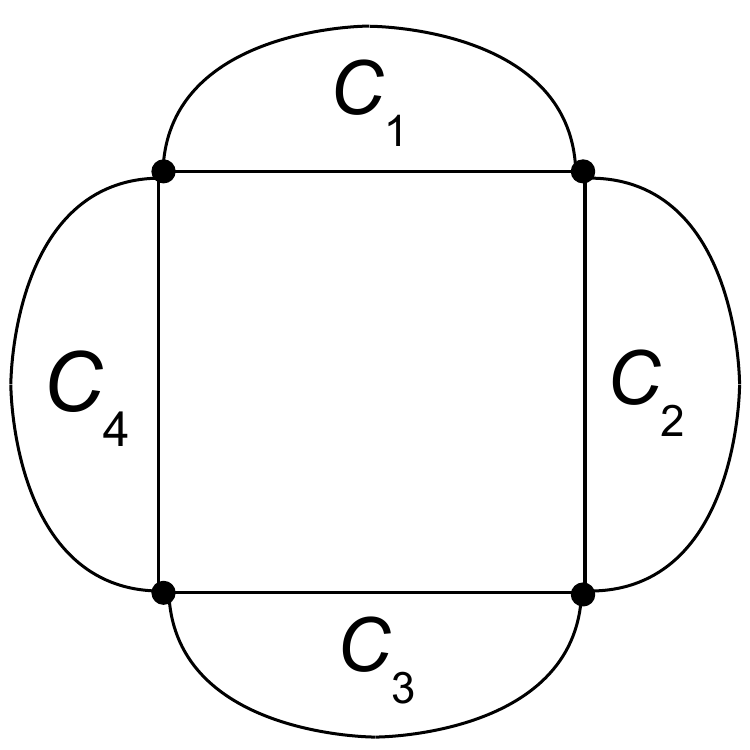}
 \caption{The graph $D_4$.
 \label {fig-D4}
 }
\end{figure}

If $D_4$ is embedded in $S^3$ in such a way that
$\lk(C_1, C_3) \ne 0$ and
$\lk(C_2, C_4) \ne 0$,
then we say it is a \dfn{\dl} $D_4$. 
If these linking numbers are both nonzero in $\zz_2$,
we say the embedded $D_4$ is \dfn{\dl\ in $\zz_2$}.
The following lemma was proved by Taniyama and Yasuhara~\cite{ty} in $\zz$
and  by Foisy~\cite {foisy-D4} in $\zz_2$.

\lemx{
\textnormal{($D_4$-Lemma,
\cite{foisy-D4},~\cite{ty})}
Every \dl\ $D_4$ (\zzt) 
contains a knot $K$ such that
$a_2(K) \ne 0$ (\zzt, respectively).
}

A \dfn{\dless} (\zzt)  of a graph is an embedding
that contains no \dl\ $D_4$ (\zzt) as a minor
(note: a minor $H$ of an embedded graph $G$
can in a natural way be given an embedding
that is ``derived'' from the given embedding of $G$).
Since any knot $K$ with $a_2(K) \ne 0$ is nontrivial,
an immediate corollary of the above lemma is:

\corx{
If a graph has no \dless\ (\zzt), then it is IK.
}

Given any graph $G$, 
our algorithm searches for a \dless\  of $G$ \zzt.
If it finds one, it reports this embedding,
which, as we'll explain later, 
may help in trying to find a knotless embedding for $G$.
If the algorithm determines that there is no such embedding,
we conclude by the above corollary that $G$ is IK.

% Using this algorithm,
% Goldberg, Mattman, and Naimi~\cite{gmn}
% have found over 1600 new IK graphs,
% most of which seem to be minor minimal IK,
% and 222 of which were ``manually'' verified in~\cite{gmn}
% to be minor minimal IK.
% 

 Foisy~\cite{foisy-ques} has asked the following question
 (which is the converse of the above corollary):

 \begin{quote}
{Question 1
\textnormal{\cite{foisy-ques}:}
 Does every graph that has a \dless\ (\zzt) have a knotless embedding?
 }
 \end{quote}

 If the answer to this question is affirmative,
 then in fact our algorithm
 (without the third shortcut described in Section ~\ref{section-shortcuts})
  \emph{always} decides
 whether or not a graph is IK.
 In~\cite{gmn}, in every instance that 
 the algorithm found a \dless\ for a graph $G$,
 whenever the authors tried to manually find a knotless embedding for $G$,
 they succeeded.
 This provides some evidence that
 the answer to Question~1 may be affirmative.

 Let's define a spatial graph $\Gamma$ to be an \dfn{$a_2$-less embedding} (\zzt)
 if for every knot $K \subset \Gamma$,
 $a_2(K) = 0$ (\zzt).
It is then natural to ask: 
 
  \begin{quote}
{Question 2:
  Does every graph that has a \dless\ (\zzt) have an $a_2$-less embedding (\zzt)?
}
\end{quote}
 
  \begin{quote}
{Question 3:
  Does every graph that has an $a_2$-less embedding (\zzt) have a knotless embedding?
 }
\end{quote}

The algorithm sets up systems of linear equations
 involving linking numbers of certain pairs of cycles in $G$
 in terms of variables representing the number of crossings between pairs of edges.
 These equations are such that
 there is a solution to at least one of the systems of equations
 if and only if $G$ has a \dless.
 The number of equations the algorithm needs to solve
 grows exponentially with the size of $G$.
 So it was necessary to find various shortcuts
 to speed up the algorithm.
 We describe the algorithm in Section~\ref{section-algorithm},
 the shortcuts in Section~\ref{section-shortcuts},
 and the proofs in Section~\ref{section-proofs}.

 Using similar techniques, 
 we have also written a Mathematica program
 that determines whether or not any given graph is intrinsically linked (IL),
 i.e., every embedding of it in $S^3$ contains a nontrivial link.
 The algorithm relies on the work of 
 Robertson, Seymour, and Thomas~\cite{rst},
 which implies that
 if a spatial graph contains a nontrivial link,
 then it contains a 2-component link with odd linking number.
 Note that, by \cite{rst}, an algorithm
 for recognizing IL graphs was already known: 
simply check whether the given graph contains as a minor
one of the seven graphs in the Petersen Family ---
which, by \cite{rs-poly}, can be done in polynomial time
(though the polynomial has very large coefficients).
And in \cite{v}, van der Holst
gives polynomial-time algorithms to
find an embedding of a graph with all linking numbers zero,
given that such an embedding exists.
%Our algorithm is not polynomial time
%(since the number of cycles can grow exponentially with the number of edges);
%but, in practice, we found it satisfactorily fast and practical
%for the graphs that we tried out:
%for graphs with 28 or fewer edges, 
%the program usually took from about a few seconds to a couple of hours to run 
%on a typical personal computer in 2009.
%More specifically, the computer program was run 
%on at least the following graphs listed in \cite{gmn} and  \cite{gmn-appendix}:
%all graphs in the $K_7$, $K_{3,3,1,1}$, and $E_9 + e$ families  (188, all together);
%the 25 parentless (``top level'') graphs and a few dozen other graphs in the $G_{9,28}$ family;
%and the graph $G_{14,25}$ plus a few more in its family.

 \section{The Algorithm}
  \label{section-algorithm}
 {
 A \dfn{quad} in a graph $G$ is an ordered 4-tuple  $\qd$
 of distinct cycles in $G$ such that
 $C_1 \cap C_3 = C_2 \cap C_4 = \emptyset$, and
 for each $i \in \{1,2,3,4\}$,
 either $C_i \cap C_{i+1}$ is nonempty and connected
 (where an index of $4+1$ is to be interpreted as 1),
 or $C_i \cap C_{i+1} = \emptyset$ and there exists a path $p_{i}$ in $G$
with one endpoint in $C_i$ and the other in $C_{i+1}$,
satisfying the following conditions:
 \begin{enumerate}
 \item
Each path $p_{i}$ is disjoint from $C_j$ for $j \not \in \{i,i+1\}$.
\item
The interior of each $p_{i}$ is disjoint from $C_i$ and $C_{i+1}$.
\item
Any two distinct paths $p_{i}$ and $p_{j}$ are disjoint.
 \end{enumerate}
 Furthermore,  for each $i$ where $C_i \cap C_{i+1}$ is nonempty and connected,
we let $p_i$ consist of an arbitrary vertex in $C_i \cap C_{i+1}$.
(We do this so that $p_i$ is defined for every $i$,
allowing for more uniformity in notation).
We say $p_1, \cdots, p_4$ are \dfn{connecting paths}
for the given quad.
 }
 \figref{example-quad-with-paths} shows an example of this. 
It is easy to see that
for any quad $\qd$ together with connecting paths $p_{i}$, 
their union $(\bigcup_i C_i) \cup (\bigcup_{i}p_{i})$ contains $D_4$ as a minor.

 \begin{figure}[ht]

 \centering
 \includegraphics[width=60mm]{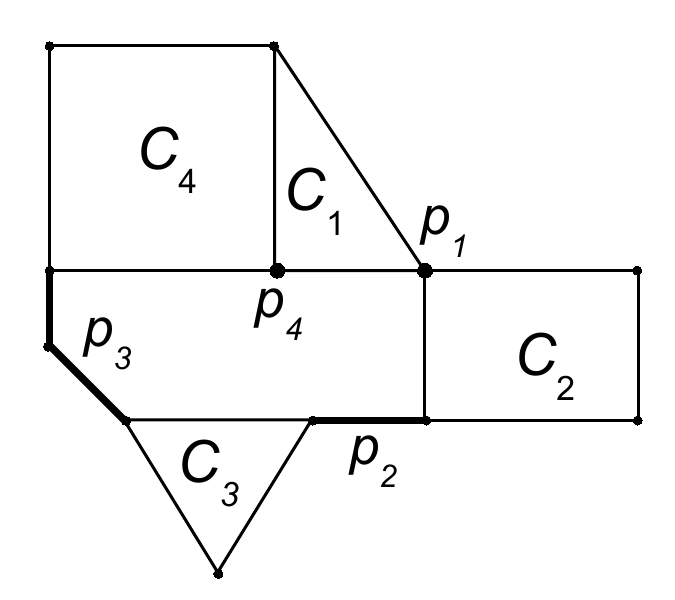}
 \caption{An example of a quad, with thickened connecting paths.
 \label {example-quad-with-paths}
 }
\end{figure}

We say a quad $\qd$
in a spatial graph $\Gamma$ 
is \dfn{\dl} (\zzt)
if
$\lk(C_1,C_3) \ne 0$ and
$\lk(C_2,C_4) \ne 0$ 
(\zzt, respectively).

It is easy to see that a spatial graph that contains a \dl\ quad
contains a \dl\ $D_4$ as a minor.
Lemma~\ref{lem-D4-quad},
 which we will prove later,
says that the converse of this is also true, i.e.,
\begin{quote}
If a spatial graph $\Gamma$ contains a \dl\ $D_4$ (\zzt) as a minor,
then it contains a \dl\ quad (\zzt, respectively).
\end{quote}
Hence an embedding $\Gamma$ of a graph $G$
is a \dless\ iff
for every quad $\qd$ in $\Gamma$,
$\lk(C_1, C_3)=0$ or $\lk(C_2,C_4)=0$.
We translate this into systems of linear equations as follows.
Let $\Gamma_0$ be a fixed arbitrary embedding of $G$.
Then any embedding $\Gamma$ of $G$ can be obtained from $\Gamma_0$
by a sequence of crossing changes together with ambient isotopy in $S^3$.
Let $e_1, \cdots, e_n$ denote the edges of $G$.
For each pair of nonadjacent edges $e_i$ and $e_j$,
let $x_{ij}$ denote the total number of crossing changes between $e_i$ and $e_j$
in the given sequence of crossing changes that takes $\Gamma_0$ to $\Gamma$.
(Crossing changes between adjacent edges are ignored since they do not affect any linking numbers.)
Note that $x_{ij}$ is the \emph{algebraic} sum of the crossing changes between $e_i$ and $e_j$,
where we give each edge an arbitrary fixed orientation, 
and associate a $+1$ to each crossing change
that changes a left-handed crossing into a right-handed one,
and a $-1$ for the opposite case (see \figref{fig-crossing-changes}).
(Alternatively, one can think of each crossing change
as adding a full positive or negative twist between the two edges.)

 \begin{figure}[ht]

 \centering
 \includegraphics[width=80mm]{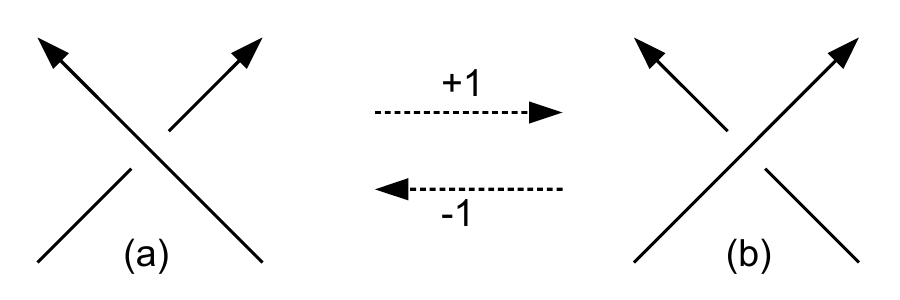}
 \caption{Changing a left-handed crossing (a) into a right-handed crossing (b) is denoted by $+1$;
 the reverse is denoted by $-1$.
 \label {fig-crossing-changes}
 }
\end{figure}

For each pair of disjoint cycles $C_i$ and $C_j$ in $G$,
denote
 the linking number of $C_i$ with $C_j$ in $\Gamma_0$ and $\Gamma$
 by $\lk(C_i,C_j,\Gamma_0)$ and  $\lk(C_i,C_j,\Gamma)$, respectively.
 Then 
 $$\lk(C_i,C_j,\Gamma) = 
 \lk(C_i,C_j,\Gamma_0) + \sum_{e_k \in C_i, \, e_l \in C_j} x_{kl}$$ 
It follows that $G$ has a \dless\ $\Gamma$ iff
for each quad $\qd$,
at least one of the following equations is satisfied:

\begin{equation}
\label{eqn:C1C3}
 \lk(C_1,C_3,\Gamma_0) + \sum_{e_k \in C_1, \, e_l \in C_3} x_{kl} = 0 
\end{equation}
or
\begin{equation}
\label{eqn:C2C4}
 \lk(C_2,C_4,\Gamma_0) + \sum_{e_k \in C_2, \, e_l \in C_4} x_{kl} = 0 
\end{equation}
where the linking numbers
 $\lk(C_i,C_{i+2},\Gamma_0)$ 
 are calculated numerically by the algorithm.
Thus, we have:

\propx{ Suppose a graph $G$ has $Q$ quads.
Choosing one of equations~(\ref {eqn:C1C3}) and (\ref {eqn:C2C4}) for each quad
gives $2^Q$ systems of linear equations,
with each system consisting of $Q$ equations,
such that
 $G$ has a \dless\ (\zzt)
iff at least one of these systems of equations has a solution (\zzt, resp.).
}
Hence, if $G$ has no \dless,
the algorithm would \emph{a priori} need to check 
every one of these $2^Q$ systems of equations to reach that conclusion.
For the graphs we worked with,
the number $Q$ of quads typically ranged 
from a few thousand to hundreds of thousands,
making $2^Q$ much too large 
to allow this algorithm to be of practical use.
Before we describe the shortcuts we used to speed up the algorithm,
we discuss two more issues:
(i)~integer solutions versus $\zz_2$ solutions, and 
(ii)~what the solutions to the equations tell us.

Solving a system of equations in $\zz$, 
as opposed to in $\zz_2$ or in $\rr$,
requires different methods and can be much slower.
So our algorithm first tries to solve each system of equations in $\zz_2$.
If there is no solution in $\zz_2$,
we know that there is no solution in $\zz$ either;
so the algorithm moves on to the ``next'' system of equations.
However, if a system of equations turns out to have a solution in $\zz_2$,
then the algorithm determines whether the system
also has a solution in $\rr$.
If it does, and this solution happens to be in $\zz$,
the algorithm outputs this integer solution
and reports that $G$ has a \dless\ in $\zz$.
If the system has no real solution
or if a real solution is found that is not in $\zz$,
the algorithm moves on to the ``next" system of equations.
In the end,
if no $\zz_2$ solutions were found for any of the systems,
the algorithm reports that the graph is IK.
If at least one $\zz_2$ solution was found but never any integer solutions
(we do not recall coming across such a case in our trials),
we only conclude that
$G$ has a \dless\ in $\zz_2$,
(although it may also have one in $\zz$).

Now let us discuss the case when the algorithm 
finds an integer solution for one of the systems of equations.
This solution tells us 
how many crossing changes to make
between each pair of nonadjacent edges
in order to go from $\Gamma_0$ to a \dless\ $\Gamma$.
These crossing changes can be implemented in different orders
and different ways up to isotopy,
thus yielding different embeddings of $G$
(\figref {fig-DifferentWaysForCrossingChange}).
And all of these embeddings will be $D_4$-less
since linking numbers depend only on 
the total numbers of crossings between nonadjacent pairs of edges.
Of course, some or all of these $D_4$-less embeddings might contain nontrivial knots.
However, every time we tried to,
we succeeded in manually finding a knotless embedding
by using these prescribed crossing changes as a ``rough guide''.
Our manual process was by no means algorithmic;
it involved making choices in how to implement each crossing change.
Making successful choices depended on both experience and luck!

 \begin{figure}[ht]

 \centering
 \includegraphics[width=100mm]{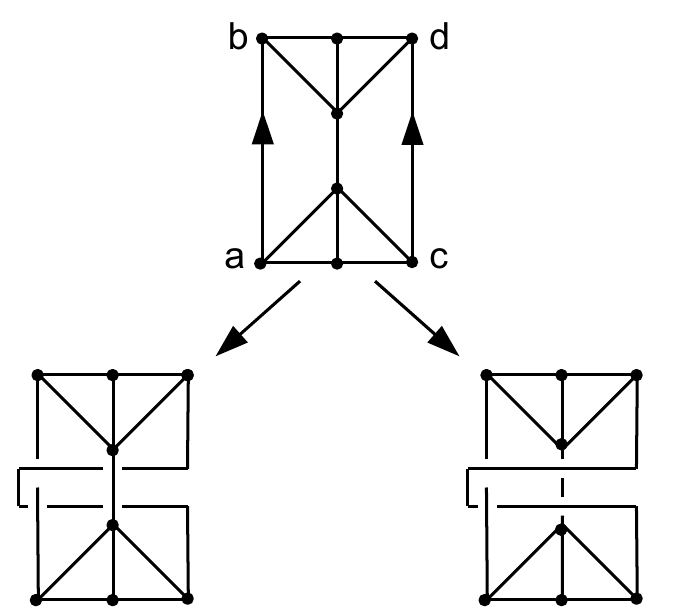}
 \caption{Two different ways (among infinitely many)
that a ``$+1$ crossing change''
(i.e., a full positive twist)
can be introduced between the oriented edges $ab$ and $cd$.
 \label {fig-DifferentWaysForCrossingChange}
 }
\end{figure}

\section{The Shortcuts}
  \label{section-shortcuts}

All of the work in this section
applies to both $\zz$ and $\zz_2$.
With this understood, we will stop repeatedly saying ``\zzt''.

Let $G$ be a graph.
For each quad $q=\qd$ in $G$,
let $q(1,3)$ and $q(2,4)$
denote the equations one gets in terms of the variables $x_{ij}$ by setting each of
$\lk(C_1, C_3) $ and  $\lk(C_2,C_4)$, respectively,
equal to zero.
Let $EQ$ be the set of all equations one obtains in this way from all quads in $G$,
i.e., 
$$EQ = \{q(i,i+2) : \;  i=1,2, \; q \mbox{ is a quad in } G \}$$
Note that each equation $A \in EQ$ may come from more than one quad.
Let's say two equations $A, B \in EQ$ are \dfn{related} 
if there is a quad $q$ such that $\{A,B\} = \{q(1,3), q(2,4)\}$.
Then, in any \dless\ of $G$, 
at least one equation in each pair of related equations must be satisfied.

One of our shortcuts is as follows.
For each equation $A \in EQ$,
we find the set $R(A)$ of all equations that are related to $A$.
(Note that by definition $A \not\in R(A)$.)
In each \dless\ of $G$,
either $A$ is satisfied or every equation in $R(A)$ is satisfied.
So if $R(A)$, considered as a system of equations, has no solution,
then $A$ must be satisfied in every \dless\ of $G$;
in this case we label $A$ as an \dfn{indispensable} equation.
We check each equation in $EQ$ for whether it is indispensable, 
and each time we find a new one,
we add it to the collection of all indispensable equations found so far
and then check whether this new collection, considered as a system of equations,
has a solution.
If the collection has no solution,
we conclude that $G$ has no \dless,
and the algorithm stops.
Even when the collection does have a solution,
this shortcut can still greatly speed up the algorithm
in either finding a \dless\ or determining that there is none:
from among the $2^Q$ systems of equations,
we work only on those systems that contain
all the indispensable equations.

Remark:
Among the graphs that we worked with,
we came across graphs $G$ of both of the following types:
(1)~$G$ has no \dless\ but its indispensable equations are satisfiable as a system;
(2)~$G$ has no \dless\ and its indispensable equations are not satisfiable as a system.

The second shortcut that we used in the algorithm is as follows.
Let $EQ = \{A_1, \cdots, A_z\}$ 
(where the number of equations, $z$,
is at most twice the number quads).
We represent each of the $2^Q$ systems of equations
by a string $S$ of 0s and 1s of length $z$,
where the $i$th digit of $S$ is 1 or 0
according to whether $A_i$ 
is or is not in the given system of equations.

Of course we are not interested in all strings of 0s and 1s of length $z$;
we say such a string is \dfn{valid} if
for each pair of related equations $A_i, A_j$,
the $i$th or $j$th digit of $S$ is 1.
The algorithm progresses through the set of all valid strings of 0s and 1s of length $z$ 
from smallest to largest in lexicographical order.
The shortcut consists of trying to skip over
collections of consecutive strings
each of which is either not valid
or corresponds to a system of equations that can be predicted,
as described below, to have no solution.

Suppose we test and determine a system of equations to have no solution. 
Often, the rows of the coefficient matrix of such a system turn out to be linearly dependent.
We find a linearly independent subset of the rows that span the rowspace.
This corresponds to a subsystem of equations that still has no solution.
Thus, any string, valid or not,
that contains 1s corresponding to this subsystem should be skipped.
This allows the algorithm to speed up greatly
by skipping over large numbers of strings.

We now describe the third shortcut we use in the algorithm.
Let $C$ be a cycle in a graph $G$, 
and suppose there exist vertices $v,w$ in $C$ such that the edge $vw$ is in $G$ but not in $C$.
Then we say the edge $vw$ is a \dfn{chord} for $C$.
If $C$ has no chords, we say it is \dfn{chordless}.
A quad is said to be chordless if all of its four cycles are chordless.
The shortcut is to only find and work with chordless quads instead of all quads.
Obviously, if a graph $G$ has a \dl\ chordless quad in every embedding,
then $G$ is IK.
However, a spatial graph that contains a \dl\ quad
may \emph{a priori} contain no \dl\ chordless quad.
So an embedding that has no \dl\ chordless quads may not be a \dless.
It appears though that such cases are rare;
in every instance that our algorithm found an embedded graph with
no \dl\ chordless quads,
whenever we tried to we succeeded in manually finding a knotless (and, hence, $D_4$-less) embedding for that graph.
%(Perhaps this is not merely by luck.
%It turns out that one can ``almost prove'' that
%if a spatial graph contains a \dl\ quad,
%then it contains a \dl\ chordless quad;
%our attempted proof failed only in one ``special'' case.)

\section{Proofs of the Lemmas}
\label{section-proofs}

Suppose $H$ is a minor of a graph $G$,
with a given sequence $S$ of edge contractions
that turn a subgraph $G'$ of $G$ into $H$.
We say a path $p'$ in $G' \subset G$ 
\dfn{corresponds} to a path $p$  in $H$
if the given sequence $S$ of edge contractions turns $p'$ into $p$.
(Note that a path may be a cycle or just one vertex.)
There may be different subgraphs $G'$ of $G$
with different sequences $S$ of edge contractions
through which one can obtain $H$ as a minor of $G$;
and correspondence between two paths certainly depends on
the choice of $G'$ and $S$.
So whenever we say $p'$ corresponds to $p$,
it is to be understood that the correspondence is relative to the given $G'$ and $S$.

When working with simple graphs,
i.e., graphs with no loops or double edges,
it is a common convention that
whenever an edge contraction results in a double edge,
the two edges in the resulting double edge
are ``merged together'' into one edge.
However, in this paper we work with $D_4$, which is a multigraph;
hence we do \emph{not} adopt this convention.
In other words, we do allow loops, and we do not merge together double edges.
An advantage of this is that
if $H$ is obtained from a subgraph of a graph $G$
by a finite sequence of edge contractions,
then clearly
\emph{to each edge in $H$ there corresponds a unique edge in $G$}.
As we  will see in Lemma~\ref {lem-unique-cycle} below
(which is used to prove Lemma~\ref {lem-D4-quad}),
the same is true for cycles (though not necessarily for all paths).

\lem{
\label{lem-unique-cycle}
Suppose $H$ is obtained from a subgraph of a graph $G$
by a given sequence of edge contractions,
where double edges are not merged together and loops are not deleted.
Then to every cycle $C$ in $H$ 
there corresponds a unique cycle $C'$ in $G$.
}

\pf{
We use induction on the number of contracted edges.
Suppose $H$ is obtained from a subgraph $G'$ of $G$
by contracting one edge $e_0$ in $G'$.
For each edge $e$ in $H$, 
let $e'$ denote the unique edge in $G'$ that corresponds to $e$.
Let $C$ be a cycle in $H$.
Let $E(C)$ denote the set of edges in $C$.
Then exactly one of
$\bigcup_{e \in E(C)} e'$
or 
$\bigcup_{e \in E(C)} e' \cup \{e_0'\}$
is a cycle in $G' \subseteq G$;
and this cycle becomes $C$ after we contract $e_0$.
The lemma now follows by induction.
}

\lem{
\label{lem-D4-quad}
If a spatial graph $\Gamma$ contains a \dl\ $D_4$ (\zzt) as a minor,
then it contains a \dl\ quad (\zzt, respectively).
}

\pf{
We use induction on the number of contracted edges.
To prove the ``base case'',
we simply note that a \dl\ $D_4$ clearly contains a \dl\ quad.

As our induction hypothesis,
assume $\Omega$ is a spatial graph 
that contains a \dl\ quad 
and is obtained from $\Gamma$
by contracting one edge $e_0 \subset \Gamma$
into a vertex $v_0 \in \Omega$.
We will show $\Gamma$ contains a \dl\ quad.

Let $\qd$ be a \dl\ quad in $\Omega$ 
with connecting paths $p_1, \cdots, p_4$.
By Lemma~\ref {lem-unique-cycle},
there exist unique cycles $C'_1, \cdots, C'_4$
in $\Gamma$ that correspond to $C_1, \cdots, C_4$.
It is easy to check that $\lk(C'_i,C'_{i+1})=\lk(C_i,C_{i+1})$ (\zzt).
We now construct connecting paths $q_1, \cdots, q_4$ for
$(C'_1, \cdots C'_4)$.
(We use $q_i$ instead of $p'_i$ to denote these paths
because some of them may not correspond to the $p_i$.)
Fix $i \in \{1,2,3,4\}$. We have two cases.

Case 1. 
$C_i \cap C_{i+1}$ is nonempty and connected.
Then one verifies easily that
$C'_i \cap C'_{i+1}$ too is nonempty and connected.
So we let $q_i$ be any vertex in $C'_i \cap C'_{i+1}$.

Case 2. 
$C_i \cap C_{i+1} = \emptyset$.
We have two subcases:

\begin{itemize}
\item[(a)]
$v_0$ is not in $p_i$, or $v_0$ is in the interior of $p_i$.
Then, by an argument similar to the proof of Lemma~ \ref{lem-unique-cycle},
there is a unique path in $\Gamma$ that corresponds to $p_i$.
We let $q_i$ be this path.

\item[(b)]
$v_0$ is in $p_i$ and $v_0$ is not in the interior of $p_i$.
Hence $v_0$ is an endpoint of $p_i$
and must lie in exactly one of $C_i$ or $C_{i+1}$.
Without loss of generality, assume $v_0 \in C_i$.
Let $p'_i$ be the path whose edges correspond to the edges of $p_i$.
If $p'_i$ has an endpoint in  $C'_i$,
we let $q_i = p'_i$;
otherwise
we let $q_i = p'_i \cup \{e_0\}$.

\end{itemize}

Now it is not difficult to see that
$(C'_1, C'_2, C'_3, C'_4)$ and $q_1, \cdots, q_4$
satisfy the definition of a \dl\ quad in $\Gamma$.
}

\paragraph*{Acknowledgment}{
The authors thank Thomas Mattman
of CSU Chico, 
and the referees, for helpful comments.
}

}\end{document}